\documentclass[12pt, leqno]{amsart}

\usepackage[usenames]{color}
\usepackage{amssymb}
\usepackage{graphicx}
\usepackage{epsf}
\usepackage{amsthm}
\usepackage{enumerate}
\usepackage[colorlinks=true,linkcolor=webgreen,filecolor=webmaroon,
  citecolor=webred]{hyperref}
\definecolor{webgreen}{rgb}{0,.5,0}
\definecolor{webbrown}{rgb}{.6,0,0}
\definecolor{webred}{rgb}{.9,.1,0}
\definecolor{webmaroon}{rgb}{0, 0.87, 0.68}
\usepackage{enumerate}
\usepackage{marvosym}
\usepackage{mathrsfs}
\usepackage{array}

\newtheorem{thm}{Theorem}

\newtheorem{lem}[thm]{Lemma}

\newtheorem{prop}[thm]{Proposition}

\theoremstyle{definition}

\numberwithin{equation}{section}

\begin{document}


\baselineskip=17pt


\title[Further Estimates with Pseudogamma functions] 
{Estimating the ratios of the Riemann Xi function 
over the pseudogamma functions}

\dedicatory{{\bf Dedicated to our Alma Mater: Wuhan University}}

\pagestyle{myheadings}

\author{Yuanyou Cheng*}
\address{\scriptsize  
Yuanyou Cheng\hfill\break
\indent 
Dept. of Mathematics and Statistics, Boston University, Boston, MA 02215 
and Dept. of Mathematics, Harvard University, Cambridge, MA 02138, USA.
\hfill\break\indent
{\tt Current Address}: Department of Mathematics, the University 
of Maryland, College Park, MD 20742, USA. \hfill}
\email{\scriptsize 
Yuanyou (Furui) Cheng \hskip -.08true cm <\hskip -0true cm  
Cheng.Y.Fred@math.umd.edu \hskip -.16true cm> } 

\author{Gongbao Li}
\address{\scriptsize  
Gongbao Li\hfill\break
\indent 
School of Mathematics and Statistics, Huazhong Normal University, 
Wuhan 430079, China. \hfill
}
\email{\scriptsize 
G. Li \hskip -.08true cm 
<\hskip -0true cm  
ligb@mail.ccnu.edu.cn 
\hskip -.16true cm> } 

\author{Juping Wang}
\address{\scriptsize 
Juping Wang\hfill\break
The School of Mathematical Sciences, Shanghai, 200433, China. 
	\hfill}
\email{\scriptsize 
J. Wang \hskip -.08true cm 
<\hskip -0true cm  jpwang@fudan.edu.cn  
\hskip -.16true cm> }

\thanks{*The authors would like to thank George E. Andrews, John T. Tate Jr., 
Roger D. Heath-Brown, Juping Wang, Andrew J. Granville, Bonnie J. Xiong, Edward 
A. Azoff, Bill Z. Yang, and Min Tang for their encouragements during the writing 
of this article. }

\begin{abstract}
A family of pseudo-Gamma functions is an essential tool, see \cite{CAPGA},
introduced recently by Y. Cheng, G. J. Fox, C. B. Pomerance, and S. W. 
Graham  and finished as in \cite{CAPGB} by Y. Cheng, C. B. Pomerance, G. J. 
Fox, and S. W. Graham for the proof of the strong density hypothesis, see 
\cite{Dh7}. These functions are reflectively symmetric about the real axis 
and also reflectively symmetric about the half line where the real part of 
the variable is equal to $\tfrac{1}{2}$. In this paper, we extend estimates 
in the proof of the strong density hypothesis towards a proof of 
the Lindel\"of hypothesis in its zero-rate form presented in \cite{Lh1}. 
Both results in \cite{Dh7} and \cite{Lh1} are used in proving 
the Riemann hypothesis in \cite{CGGP}.   
\end{abstract}
 
\date{Drafted on December 31, 2012.  Revised on July 19, 2018. 
Submitted, February 9, 2022. \hfil\\ 
\\
{\tt Ethical Statement}: {\em This article is 
an original work. It has been submitted only to the Journal of Analysis for 
publication. All authors are solemnly holding ethical standard by the journal 
concerning this work. } \\
{\bf Conflict of Interest}: The authors declare that they have no conflict of 
interest. \hfil\\ }



\maketitle
\vskip -1.4true cm 

\markleft{Yuanyou Cheng, Gongbao Li, and Juping Wang}

\thispagestyle{empty}

\section{Introduction}
\label{sec: intro}

\noindent 
We introduced a family of pseudogamma functions in \cite{CAPGA} and \cite{CAPGB}, 
mainly for the purpose of proving the density hypothesis in \cite{Dh7}, with 
more required estimates involving this family of pseudogamma functions 
in \cite{CAPGA} and \cite{CAPGB}.

We recall that the density hypothesis gives an estimate on the number 
$N(\lambda, T)$ of zeros of the Riemann zeta function $\zeta(s)$, 
$s =\sigma +i t$, in the region such that $\sigma >\lambda$,  $0 <t <T$, 
for any $\tfrac{1}{2}< \lambda <1$ and $T \ge T\sb{0}$, with $T\sb{0}$ 
sufficiently large. 

In \cite{Lh1} we prove that the Lindel\"of hypothesis, which gives a bound 
on $\zeta\bigl( \tfrac{1}{2} +i t\bigr)$ for $|t|$ large, by verifying its 
equivalent form in the zero-rate form with $N(\lambda, T +d) -N(\lambda, 
T -d)$ (the equivalence is well-known in the literature), for $1 \le d 
\le 5/4$. 

The pseudogamma function is defined in \cite{CAPGA} and is given by 
\begin{equation}
\label{eq: nabladefi}
\nabla(s) =\omega\Bigl( \tfrac{ W\sb{2} - 1/2 }{ W\sb{1} -1/2 } \Bigr)\sp{\!\!q}\, 
\biggl[ \prod\sb{k=1}\sp{ 2\sp{ K+1} } {\mathscr P}(k; s) \biggr]
\sp{\frac{q}{ 2\sp{K+1} }  }
\end{equation} 
where 
\begin{equation}
\label{eq: prenabla}
{\mathscr P}(k; s) =\dfrac{\, (s -\tfrac{1}{2})  
 -e\sp{\frac{i\,k \pi}{ 2\sp{K} } } (W\sb{1}-\tfrac{1}{2}) \,} 
{\, (s -\tfrac{1}{2})  -e\sp{\frac{i\,k \pi}{ 2\sp{K} } }
(W\sb{2} -\tfrac{1}{2}) \,},
\end{equation}
with $W\sb{1} =3 R +R\sp{1/4} +\tfrac{1}{2}$, $W\sb{2} =3 R +\tfrac{1}{2}$, 
${\grave\gamma} =0.3674$, where $\omega=\xi(1/2)$ with $\xi(1/2) >0.497$, 
$\Omega=\tfrac{47.545}{\xi(1/2)}$, and 
\begin{equation*}
\begin{split}
q &=\dfrac{(R +10)\log (R/2) +4\log\Omega }{ 4\ {\grave\gamma}\, R\sp{1/4}}, \\
	K &=\biggl\lfloor \dfrac{\log \bigl[ \tfrac{10}{27}
 		( 15 +  \tfrac{4}{R\sp{3/4}} \bigr) R\sp{1/4} +\tfrac{2}{3 R\sp{3/4}}  
 			\bigr] +2\log R}{\log 2} \biggr\rfloor,\\
\end{split}
\end{equation*}
as in Theorem 1 of \cite{CAPGA}, and $R\ge R\sb{0} =4891999112053$. In fact 
the proofs in \cite{Dh7} use $2 T -1 < R\le 2 T+1$, $T \ge T\sb{0}$ with $T\sb{0} 
=2445999556027$ and the aobe $R\ge R\sb{0}:= 2 T\sb{0} -1$.  It is obvious that 
$\nabla(s)$ is reflectively symmetric about the real axis as well as about 
the line $\sigma =\tfrac{1}{2}$; and $\nabla(\tfrac{1}{2})=\xi(1/2)$, from 
the above definition. 

For the proofs in \cite{Dh7} and Theorem 1 in \cite{Lh1} we use in an essential 
way the new method involving the above pseudo-Gamma function $\nabla(s)$, which 
amounts in exploiting the double symmetry property of the Riemann $\xi$-function 
and the function $\nabla(s)$, to counteract the ``unwanted'' effect caused by 
the magnitude of the classical Euler Gamma function entering the estimates of 
the Riemann $\xi$-function and thus also on the Riemann zeta function $\zeta(s)$. 
Our proof in \cite{Dh7} extend methods of Backlund's the proof the Riemann-von 
Mangoldt Theorem, with the essential new method using the function $\nabla(s)$. 

For these estimates, we also use the fact that our pseudo-Gamma function, whose 
definition depends primarily on parameter $R >0$, is free of zeros and poles 
in the regions of interest. Let us also mention that in our proofs we also make 
use of the Euler product representation of the functions $\xi(s)$ and $\zeta(s)$, 
through the logarithmic derivatives of these functions. 

One may see from the following two lemmas that the pseudogamma function $\nabla(s)$ 
has about the same ``size'' as the Riemann xi funciton $\xi(s)$ in the concerned 
region, specified in Lemma 2 \ref{lem: l310} below.  

First, we quote Lemma 2 from \cite{CAGP}: 

\setcounter{thm}{0}

\begin{lem}
\label{lem: l310}
We have 
\begin{equation}
\label{eq: xitypeupd}
|\xi(s) | \le 47.545 \bigl( \tfrac{R}{2} \bigr)\sp{1.000315(R/4 +5/2)},  
\end{equation}
on the circle $|s -u| =R$, $\tfrac{1}{2} <u \le 2$, for $R\ge R\sb{0}$, with 
$R\sb{0}$ as defined above.  
\end{lem}

In Section \ref{sec: pfthm2}, after recalling some relevant results of 
\cite{CAPGB}, we state and prove all results in this article, except 
the proof of Lemma \ref{lem: BC12updc}, which is postponed to Section 
\ref{sec: sec3}.  

In Section \ref{sec: sec3}, we give the proof of Lemma \ref{lem: BC12upd} 
in five subsections for clarity, as the computations involved are 
complicated.

\section{\normalsize The main results}
\label{sec: pfthm2}

In the following lemma, we quote Theorem 2 from \cite{CAPGA}. 

\begin{lem}
\label{lem: forLh}
For the pseudogamma function defined in \eqref{eq: nabladefi}, we have 
\begin{equation}
\label{eq: Omegaalpha1}
47.545\, \bigl( \tfrac{R}{2} \bigr)\sp{\frac{R}{4} +\frac{5}{2} }\,  
<|\nabla(s)| <47.545\, \bigl( \tfrac{R}{2} \bigr)\sp
{1.000315 \bigl( \frac{R}{4} +\frac{5}{2}\bigr) }\,,
\end{equation}
on the circle $| s -\tfrac{1}{2} | =R$, $R \ge R\sb{0}$.
\end{lem}

We need the following result from Lemma 5 in \cite{CAGP}. 

\setcounter{thm}{0}

\begin{prop}
\label{prop: mainlem}
Let $\nabla(s)$ be the pseudo-Gamma function defined in \eqref{eq: nabladefi} 
with the choice of constants presented after \eqref{eq: prenabla} and let 
$R \ge 2\,T -1$ with $T \ge T\sb{0}$, where $T\sb{0}$ is defined as above. 
Then, $\nabla(u)$ is a monotonously increasing function of $u\in(1/2, 2]$ 
and 
\begin{equation}
\label{eq: intact}
\omega\le \nabla(u) \le R\sp{\frac{1}{4.408\, R}} \le \bigl( 1 +10\sp{-11} \bigr)
\omega, 
\end{equation}
for $\tfrac{1}{2} \le u \le 2$. \qed
\end{prop}

The following lemma is the content of Lemma 6 and  Lemma 4 from \cite{CAGP}. 

\begin{prop}
\label{prop: l31} 
Let $T\ge T\sb{0}$ and $2 T -1 < R <2 T +1$ with $T\sb{0}$ defined in the second 
paragraph in Section \ref{sec: intro}. Suppose that there are no zeros for 
the Riemann xi-function $\xi(s)$ on the circle $\bigl| s -u \bigr| =R$, with 
$\tfrac{1}{2}< u \le 2$. Then the function ${\mathbf B}(s):=\tfrac{\xi(s)}
{\nabla(s)}$ is analytic inside any circle $|s-u| =R +\varepsilon$, with 
$\tfrac{1}{2}< u\le 2$ and sufficiently small $\varepsilon >0$, and satisfies 
the following upper bound  
\begin{equation}
\label{eq: nablaBs} 
{\mathbf B}(s) \trianglelefteq b\sb{0},
\end{equation}
with $b\sb{0} =1.0132$, on the circle $\bigl|s - u\bigr|= R$, with $\nabla(s)$ 
in \eqref{eq: nabladefi} with \eqref{eq: prenabla} and the choice of parameters 
specified after \eqref{eq: prenabla}. Furthermore, we have 
\begin{equation}
\label{eq: dot497}
\xi(u) >\xi\bigl( \tfrac{1}{2} \bigr) >0.497,
\end{equation} 
for $\tfrac{1}{2} <u \le 2$. \qed
\end{prop}

We also need another lemma, stated below as Lemma \ref{lem: BC12upd}. Let $\eta$ 
be defined as in (3.1) with $X\sb{0}$ and $Y\sb{0}$ defined in the paragraph 
after (3.1) in \cite{Dh7} and  $X =X\sb{0} +\tfrac{\eta}{2}$, $Y=Y\sb{0}
+\tfrac{\eta}{4}$, and $Y\sb{1} =Y\sb{0} +\tfrac{3\eta}{4}$ where $T -\tfrac{\eta}
{4} <Y <Y\sb{1} <T +\tfrac{3\eta}{4}$. We recall the definition of ${\mathbf B}(s)$ 
in Lemma \ref{prop: l31} and set
\begin{equation}
\label{eq: BC12defi}
 {\mathbf C}(s) =\tfrac{ \nabla( 3/2 + s) }{ \nabla( X -1/2 +s)}, \quad
 {\mathbf D}(s) =\tfrac{\nabla(s +i\,Y\sb{1})}{\nabla(s +i\,Y)}, 
\end{equation}
where the variables $X$, $Y$, and $Y\sb{1}$ introduced in the above are the same 
as those used in Section 3 in \cite{Dh7}.

\begin{lem} 
\label{lem: BC12upd}
Let $T\ge T\sb{0}$ with $T\sb{0}$ defined in Section \ref{sec: intro}, 
$ 2 T -1 < R \le 2 T +1$, $\omega =\xi(1/2)$, $\Omega =\tfrac{47.545}
{\xi(1/2)}$, and $\alpha=10$. Also, let $\sigma\sb{0}=\tfrac{5}{4}$, 
$T-\tfrac{3}{2} <t\sb{0} \le T +\tfrac{3}{2}$, and $s\sb{0} =\sigma\sb{0} 
+i\, t\sb{0}$. For the function ${\mathbf B}(s)$ defined in Lemma 
{\rm \ref{prop: l31}}, we have 
\begin{equation}
\label{eq: nablaBs1}
\tfrac{ {\mathbf B}(s)}{{\mathbf B}( 1/ 2)} \trianglelefteq 1.014,
\quad \tfrac{ {\mathbf B}(s)}{{\mathbf B}( \sigma\sb{0} )} 
 \trianglelefteq 1.014, 
\end{equation}
on the circle $| s -\tfrac{1}{2}| =R$ and $| s -\sigma\sb{0}| =R$, 
respectively.  
\end{lem}

\begin{proof}
For the proof of the first estimate in \eqref{eq: nablaBs1}, it suffices to 
recall that $|{\mathbf B}(s) |\le b\sb{0}$ from \eqref{eq: nablaBs} with 
the value of $b\sb{0}$ defined right after the inequality, $\xi(\tfrac{1}{2}) 
>0.497$ from \eqref{eq: dot497}, and $\nabla(\tfrac{1}{2}) =\xi(1/2)$, by 
the very construction of $\nabla(s)$. Therefore, 
$\bigl| \tfrac{{\mathbf B}(s) }{ {\mathbf B}(1 /2)} \bigr| <1.0132$ as 
${\mathbf B}(1/2) =1$, by the definition of ${\mathbf B}(s)$, which proves 
the first estimate in \eqref{eq: nablaBs1}. 

In order to prove the second estimate in \eqref{eq: nablaBs1}, one uses 
\eqref{eq: nablaBs} and \eqref{eq: dot497} with $R\le 2 T +1$, observing that 
$\tfrac{{\mathbf B}(s)}{{\mathbf B}(\sigma\sb{0})} =\tfrac{\xi(s)\, 
\nabla(\sigma\sb{0}) }{ \xi(\sigma\sb{0})\, \nabla(s) } \trianglelefteq 
\tfrac{(1+10\sp{-11} )\omega}{0.497} {\mathbf B}(s) <1.014$ by 
\eqref{eq: dot497},  and $\nabla(\sigma\sb{0}) <\bigl( 1+10\sp{-11} \bigr)
\omega$ from \eqref{eq: intact}, and \eqref{eq: nablaBs}.
\end{proof}

\begin{lem} 
\label{lem: BC12updb}
Let $\sigma\sb{0}=\tfrac{5}{4}$, $T-\tfrac{3}{2} <t\sb{0} \le T +\tfrac{3}{2}$, 
and $s\sb{0} =\sigma\sb{0} +i\, t\sb{0}$. For the functions ${\mathbf C}(s)$ 
and ${\mathbf D}(s)$, we have 
\begin{equation}
\label{eq: nablaCsDs}
\tfrac{{\mathbf C}(s)}{{\mathbf C}(s\sb{0})} 
 \trianglelefteq T\sp\frac{11}{T\sp{1/4}}, \quad
 \tfrac{{\mathbf D}(s)}{{\mathbf D}(s\sb{0})} 
 \trianglelefteq T\sp{\frac{11}{T\sp{1/4} } }, 
\end{equation}
on the circle $| s -s\sb{0} | =1$. 
\end{lem}

The proof of this lemma is a little more complicated; we put it in the next section. 

\begin{lem} 
\label{lem: BC12updc}
Let $\sigma\sb{0}=\tfrac{5}{4}$, 
$T-\tfrac{3}{2} <t\sb{0} \le T +\tfrac{3}{2}$, and $s\sb{0} =\sigma\sb{0} 
+i\, t\sb{0}$. We have 
\begin{equation}
\label{eq: nablaBs2}
\tfrac{ {\mathbf B}(s -i(Y\sb{1} -Y) )}{{\mathbf B}( 1/ 2)} 
 \trianglelefteq 1.015, \quad \tfrac{ {\mathbf B}
 (s -i(Y\sb{2} -Y\sb{1}) )} {{\mathbf B}( \sigma\sb{0} )} 
 \trianglelefteq 1.015, 
\end{equation}
with $T -\tfrac{1}{4} <Y\sb{1} <Y\sb{2} \le T +\tfrac{1}{4}$ {\rm(}as 
in Proposition 1 and Proposition 2 in \cite{Dh7}{\rm)}, 
on the circle $| s -\tfrac{1}{2}| =R$ and $| s -\sigma\sb{0}| =R$, 
respectively, and 
\begin{equation}
\label{eq: nablaCsDs2}
\tfrac{{\mathbf C}(s -i(Y\sb{2} -Y\sb{1}) )}{{\mathbf C}(s\sb{0})} 
 \trianglelefteq T\sp\frac{11.001}{T\sp{1/4}}, \quad
 \tfrac{{\mathbf D}(s -i(Y\sb{2} -Y\sb{1}) )}{{\mathbf D}(s\sb{0})} 
 \trianglelefteq T\sp{\frac{11.001}{T\sp{1/4} } }, 
\end{equation}
on the circle $| s -s\sb{0} | =1$. 
\end{lem}

\begin{proof}
For the estimate in \eqref{eq: nablaBs2} in Lemma \ref{lem: BC12upd}, we recall 
the definition of $\eta$ defined in (3.1) in \cite{Dh7} and notice that we may 
let $\eta$ be sufficiently small. We prove those two inequalities 
in \eqref{eq: nablaBs2} similarly; for example, we only give the details 
for the proof of the first one in \eqref{eq: nablaBs2}.  

We notice that $\tfrac{{\mathbf B}(s -i( Y\sb{1} -Y))}{{\mathbf B}(s)} 
\trianglelefteq 1+\epsilon$, for arbitrary $\epsilon >0$ if only we allow $\eta$ 
to be as small as necessary with respect to the value of $\epsilon$. Let us say 
that we use $\epsilon \le 1\times 10\sp{-9}$. It then follows that 
$\tfrac{{\mathbf B}(s -i( Y\sb{1} -Y))}{{\mathbf B}(1/2)} 
=\tfrac{{\mathbf B}(s -i( Y\sb{1} -Y))}{{\mathbf B}(s)}\, 
\tfrac{{\mathbf B}(s)}{{\mathbf B}(1/2) } \trianglelefteq 1.014\times (1 +10\sp{-9}) 
\le 1.015$, where we used the first inequality in \eqref{eq: nablaBs1}.  

Similarly, we can prove the last part \eqref{eq: nablaCsDs2} of Lemma 
\ref{lem: BC12upd}, using the second estimates in \eqref{eq: nablaCsDs} in 
Lemma \ref{lem: BC12updb}. This finishes the proof of Lemma \ref{lem: BC12upd}, 
once we prove Lemma \ref{lem: BC12updb}.
\end{proof}

\section{Proof of Lemma \ref{lem: BC12updb} }
\label{sec: sec3}

We now give a proof of Lemma \ref{lem: BC12updb}. To prove the first inequality  
in \eqref{eq: nablaCsDs}, with the estimates on the circle and at the center 
$s\sb{0}$ of $|s -s\sb{0} | =1$, we divide the proof into four subsections. 
We prove the second inequality in \eqref{eq: nablaCsDs} in the last subsection.  

\subsection
{Set up for introducing the logarithm of the quantities in \ref{eq: nablaCsDs}}. 
\label{subsec: 31}

 From the definition of ${\mathbf C}(s)$ in \eqref{eq: BC12defi} and of $\nabla(s)$ 
 in \eqref{eq: nabladefi}, we have 
\begin{equation}
\label{eq: overandover}
\tfrac{{\mathbf C}(s)}{{\mathbf C}(s\sb{0})} 
=\tfrac{\nabla(3/2+s)}{\nabla(3/2 +s\sb{0})}\big/ \tfrac{\nabla(X -1/2 +s)}
{\nabla(X -1/2 +s\sb{0})} .
\end{equation}
Denote 
\begin{equation}
\label{eq: s1s}
s\sb{1} =s, \quad\text{on the circle}\quad |s -s\sb{0}| =1,
\end{equation}
for brevity, without repeating several expressions for both $s\sb{0}$ and $s\sb{1}$. 
We should keep in mind that the subscript $1$ stays for the point on the circle and 
the subscript $0$ is for the point at its center. Recalling the expression 
\eqref{eq: overandover} for $\tfrac{{\mathbf C}(s)}{ {\mathbf C}(s\sb{0}) }$, 
the variable inside $\nabla$ for the numerator $\nabla(\tfrac{3}{2} +s)$ is equal 
to $\tfrac{3}{2} +s\sb{l} -\tfrac{1}{2} =s\sb{l} +1$ for $l=0$ and $1$ and that 
for the denominator $\nabla(X -\tfrac{1}{2} +s)$ is equal to $X -\tfrac{1}{2} 
+s\sb{l} -\tfrac{1}{2} =s\sb{l} -(1 -X)$, respectively. We write for $l=0$ and $1$ 
\begin{equation}
\label{eq: calculatelater}
\begin{split}
s\sb{l} +1 &=\sigma\sb{l} +1 +i\,t\sb{l} 
 ={\hat r}\sb{l}\, e\sp{i {\hat\theta}\sb{l}}, \\
s\sb{l} -(1 -X) &=\sigma\sb{l} -(1 -X) +i\,t\sb{l} 
 ={\check r}\sb{l}\, e\sp{i {\check\theta}\sb{l}}, \\
\end{split}
\end{equation} 
with the `hat' on the top of variables for the numerators and the `check' on 
the top of variables for the denominators, in both the rectangular and polar 
coordinate systems, respectively. Note that ${\hat r}\sb{l} \ge 0$, 
${\check r}\sb{l} \ge 0$, $0\le {\hat\theta}\sb{l} <2\pi$, and $0 \le 
{\check\theta}\sb{l} <2\pi$, $l =0$ and $1$. We have, using $\sigma\sb{1} 
=\sigma$ and $t\sb{1} =t$, from \eqref{eq: s1s}:
\begin{equation}
\label{eq: eight}
\begin{split}
&{\hat r}\sb{0} =|s\sb{0} +1|, 
  \hskip 3.2true cm
  {\hat\theta}\sb{0} =\tfrac{\pi}{2} -\arctan\tfrac{9}{4t\sb{0}};  \\
&{\hat r}\sb{1} =| s\sb{1} +1 |, 
 \hskip 3.2true cm 
 {\hat\theta}\sb{1} =\tfrac{\pi}{2} -\arctan\tfrac{\sigma +1}{t};  \\
&{\check r}\sb{0} =|s\sb{0} -( 1 -X) |, 
 \hskip 2true cm
 {\check\theta}\sb{0} =\tfrac{\pi}{2} -\arctan\tfrac{1/4 +X}{t\sb{0}}; \\
&{\check r}\sb{1} =|s\sb{1} -( 1 -X) |, 
 \hskip 2true cm 
 {\check\theta}\sb{1} =\tfrac{\pi}{2} -\arctan\tfrac{\sigma -(1 -X)}{t}, \\
\end{split}
\end{equation}
recalling that $s\sb{0} =\sigma\sb{0} +i\,t\sb{0}$ from the statement of 
Lemma \ref{lem: BC12upd} and using $\arctan(x) = \tfrac{\pi}{2} 
-\arctan(\tfrac{1}{x})$ for $x >0$. To reduce the estimates in \eqref{eq: 3stuffs},
\eqref{eq: 3stuffs1} and \eqref{eq: 3stuffsg} below, it is useful to give 
the estimates for ${\hat r}\sb{l}$, ${\check r}\sb{l}$, ${\hat\theta}\sb{l}$, and 
${\check\theta}\sb{l}$, for $l =0$, $1$ in what follows.

We omit the detailed computation for the expressions above for the $G$'s and $g$'s, 
one may refer to (3.2) in \cite{CAPGA} for the details or figure it out by noting 
$|r\, e\sp{i \theta} -S\, e\sp{i \phi}|\sp{2} =r\sp{2} +S\sp{2} -2 r S \cos(\theta
-\phi)$ for any $r \in{\mathbf R}$, $S\in{\mathbf R}$, $\theta \in{\mathbf R}$
and $\phi\in{\mathbf R}$, keeping in mind that we have only used the notations 
$g$'s there corresponding to $G$'s and $g$'s here.

We set
\begin{equation}
\label{eq: ggGG1210}
\begin{split}
G\sb{j}({\hat\theta}\sb{l}, {\hat r}\sb{l}; k) 
&=|(s\sb{l} +1) -\exp(i k\pi) (W\sb{j} -1/2)|\sp{2} \\
&=| {\hat r}\sb{l} \cos{\hat\theta} +i {\hat r}\sb{l}\sin{\hat\theta}\sb{l} 
 -( \cos\tfrac{k\pi}{2\sp{K}} +i \sin\tfrac{k\pi}{2\sp{K}} )
 (W\sb{j} -\tfrac{1}{2})|\sp{2} \\
& = (W\sb{j} -\tfrac{1}{2})\sp{2} +{\hat r}\sb{l}\sp{2} -2(W\sb{j} -\tfrac{1}{2}) 
 {\hat r}\sb{l} \cos({\hat\theta}\sb{l} -\tfrac{k\pi}{2\sp{K} } ), \\
g\sb{j}(l) 
&=|[s\sb{l} -(1 -X)] -\exp(i k\pi) (W\sb{j} -1/2)|\sp{2} \\
& = (W\sb{j} -\tfrac{1}{2})\sp{2} +{\check r}\sb{l}\sp{2} -2(W\sb{j} -\tfrac{1}{2}) 
 {\check r}\sb{l} \cos({\check\theta}\sb{l} -\tfrac{k\pi }{ 2\sp{K} } ),  \\
\end{split}
\end{equation}
for $j =1$, $2$ and $l =0$, $1$, with $l =1$ corresponding to $s\sb{1} =s$ 
on the circle $|s -s\sb{0}| =1$ and $l =0$ corresponding to the center $s\sb{0}$. 

We now use the logarithm to get, using the definitions of 
${\mathbf C}(s)$ and $\nabla(s)$, 
\begin{equation}
\label{eq: logover}
\begin{split}
& \log\bigl| \tfrac{{\mathbf C}(s)}{{\mathbf C}(s\sb{0})} \bigr|
 =\bigl( \log|\nabla(3/2 +s)| -\log|\nabla(3/2 +s\sb{0})| \bigr)  \\
&\hskip 0.4true cm -\bigl( \log|\nabla(X -1/2 +s)| 
 -\log|\nabla(X -1/2 +s\sb{0})| \bigr)  \\
& =\dfrac{q}{2\sp{K+2} } \sum\sb{k =1}\sp{2\sp{K+1} } \biggl(
 \biggl[ \log\dfrac{ G\sb{1}({\hat\theta}\sb{1}, {\hat r}\sb{1}; k) }
 { G\sb{2}({\hat\theta}\sb{1}, {\hat r}\sb{1}; k) } 
 -\log\dfrac{ G\sb{1}({\hat\theta}\sb{0}, {\hat r}\sb{0}; k) }
 { G\sb{2}({\hat\theta}\sb{0}, {\hat r}\sb{0}; k) }  \biggr] \\
&\hskip 1.6true cm
  -\biggl[ \log\dfrac{ g\sb{1}({\check\theta}\sb{1}, {\check r}\sb{1}; k) }
 { g\sb{2}({\check\theta}\sb{1}, {\check r}\sb{1}; k) } 
 -\log\dfrac{ g\sb{1}({\check\theta}\sb{0}, {\check r}\sb{0}; k) }
 { g\sb{2}({\check\theta}\sb{0}, {\check r}\sb{0}; k) } \biggr] \biggr), \\
\end{split}
\end{equation}
where the $G$'s are corresponding to the numerators with the `hat' variables
and $g$'s are corresponding to the denominator with the `check' variables.

\subsection{Upper bound of the logarithm}. 
Keeping the last remarks on $G$ for numerators and $g$ for denominators 
in mind, we may write conveniently 
\begin{equation}
\label{eq: 2eqns}
\begin{split}
G\sb{j}(l) = G\sb{j} ({\hat\theta}\sb{l}, {\hat r}\sb{l}; k), \quad
\text{and}\quad g\sb{j}(l) = g\sb{j} ({\check\theta}\sb{l}, 
 {\check r}\sb{l}; k), 
\end{split}
\end{equation}
for $j =1$, $2$ and $l =0$, $1$ from now on, with the index variable 
$k$ also omitted without confusion. We notice that $\log(1+x) 
=x-\tfrac{x\sp{2}}{2} +\tfrac{x\sp{3}}{3} -\tfrac{x\sp{4}}{4}+\ldots$,
as the series is convergent for $0 <x <1$. Therefore, denoting 
\begin{equation}
\label{eq: etadefi}
\eta( x) =\tfrac{\log(1+x)}{x}, 
\end{equation}
we have $1 -\tfrac{x}{2} <\eta(x) <1$. From (3.3.7) in \cite{Dh7} with 
both $G$ and $g$ represented by $g$ there and the remark in the paragraph 
between (3.3.7) and (3.3.9) in \cite{Dh7}, we have 
\begin{equation}
\label{eq: fromCA1}
\begin{split}
G\sb{2}(l) < &G\sb{1}(l) < 2 G\sb{2}(l) , \\
g\sb{2}(l) < &g\sb{1}(l) < 2 g\sb{2}(l) , \\
\end{split}
\end{equation}
For $l=0$, $1$, we let 
\begin{equation}
\label{eq: g0G0}
\begin{split}
G\sb{0}(l) := &G\sb{1}(l) -G\sb{2}(l) . \\
g\sb{0}(l) := &g\sb{1}(l) -g\sb{2}(l) . \\
\end{split}
\end{equation}
We also have 
\begin{equation}
\label{eq: two01}
0 <\tfrac{G\sb{0}(l)}{G\sb{2}(l)} <1, \quad\text{and}\quad
 0 <\tfrac{g\sb{0}(l)}{g\sb{2}(l)} <1,
\end{equation}
for $l =0$, $1$. We shall apply the above inequality for $\eta(x)$ in 
\eqref{eq: etadefi}, with $x =\tfrac{G\sb{0}(l)}{G\sb{2}(l)}$ and 
$\tfrac{g\sb{0}(l)}{g\sb{2}(l)}$, respectively. Moreover, we have, 
with $\eta\sb{x} =\eta(x)$  for $x =\tfrac{h\sb{0}}{h\sb{2}}$ with 
$h\sb{0}$ being either $G\sb{0}(l)$ or $g\sb{0}(l)$ and $h\sb{2}$ being 
either $G\sb{2}(l)$ or $g\sb{2}(l)$,
\begin{equation}
\label{eq: only10forl}
\begin{split}
\eta\sb{G}\ \tfrac{ G\sb{0}(l) } {G\sb{2}(l)}  & \le\log\tfrac{G\sb{1}(l) } 
	{G\sb{2}(l)} \le \tfrac{ G\sb{0}(l) }{ G\sb{2}(l) }, \\
\eta\sb{g}\ \tfrac{ g\sb{0}(l) }{g\sb{2}(l)} & \le\log\tfrac{g\sb{1}(l) }
	{g\sb{2}(l)} \le \tfrac{ g\sb{0}(l) }{ g\sb{2}(l) }, \\
\end{split} 
\end{equation}
for $l =0$, $1$, and at the same time
\begin{equation}
\label{eq: usednow}
\begin{split}
1- \tfrac{ G\sb{0}(l)}{2\, G\sb{2}(l)} <\eta\sb{G} & <1, \\
1- \tfrac{ g\sb{0}(l)}{2\, g\sb{2}(l)} <\eta\sb{g} & <1. \\ 
\end{split} 
\end{equation}
It follows from 
\eqref{eq: logover} that 
\begin{equation}
\label{eq: lowerupper}
\begin{split}
\log\Bigl| \tfrac{{\mathbf C}(s)}{{\mathbf C}(s\sb{0})} \Bigr| 
 &\le \tfrac{q}{2\sp{K+2} } \Biggl[\ \sum\sb{k =1}\sp{2\sp{K+1} } \Bigl(  
  \tfrac{G\sb{0}(1) }{ G\sb{2}(1)} -\eta\sb{G}\, 
   \tfrac{ G\sb{0}(0) }{ G\sb{2}(0) } \Bigr) 
 -\sum\sb{k =1}\sp{2\sp{K+1} } \Bigl( \eta\sb{g}\, \tfrac{g\sb{0}(1) }
 { g\sb{2}(1)} -\tfrac{ g\sb{0}(0) }{ g\sb{2}(0) } \Bigr) \Biggr]. \\
\end{split}
\end{equation}

To estimate the upper bound for these two sums over the set of $k =1$, $2$, 
$\ldots$, $2\sp{K+1}$, we transform each summand in the expression on the 
right hand side of \eqref{eq: lowerupper} by the following algebraic identity. 
For all real valued variables $u$, $v$, $U$, $V$, $\eta\sb{G}$, and $\eta\sb{g}$, 
we have 
\begin{equation}
\begin{split}
\label{eq: etaGg}
\tfrac{u}{U} -\eta\sb{G}\,\tfrac{v}{V} &=\tfrac{u}{U} -\tfrac{v}{V} 
 +(1 -\eta\sb{G}) \tfrac{v}{V}, \\
\eta\sb{g}\,\tfrac{u}{U} -\tfrac{v}{V}
 &=\tfrac{u}{U} -\tfrac{v}{V} -(1 -\eta\sb{g}) \tfrac{u}{U}, \\
\end{split}
\end{equation}
\begin{equation}
\label{eq: algebraicI}
\tfrac{u}{U} -\tfrac{v}{V} = \tfrac{u V -U v}{UV} 
 =\tfrac{(u-v) V -(U -V)v}{UV}.
\end{equation}
Denote
\begin{equation}
\label{eq: h1h2}
\begin{alignedat}{2}
H\sb{1} &=G\sb{0}(1) -G\sb{0}(0), \ \quad& H\sb{2} &=G\sb{2}(1) -G\sb{2}(0), \\
h\sb{1} &=g\sb{0}(1) -g\sb{0}(0), \      & h\sb{2} &=g\sb{2}(1) -g\sb{2}(0), \\
\end{alignedat}
\end{equation}
using the definitions in \eqref{eq: 2eqns}. It follows that 
\begin{equation}
\label{eq: twoidentities}
\begin{split}
\tfrac{(G\sb{0}(1) -G\sb{0}(0) ) G\sb{2}(0) -(G\sb{2}(1) -G\sb{2}(0)) 
 G\sb{0}(0)} {G\sb{2}(1) G\sb{2}(0)} &= \tfrac{H\sb{1} G\sb{2}(0) 
 -H\sb{2} G\sb{0}(0) }{ G\sb{2}(1) G\sb{2}(0) }, \\
\tfrac{ (g\sb{0}(1) -g\sb{0}(0) ) g\sb{2}(0) -(g\sb{2}(1) -g\sb{2}(0) ) 
 g\sb{0}(0)}{g\sb{2}(1) g\sb{2}(0)} &= \tfrac{h\sb{1} g\sb{2}(0)
 -h\sb{2} g\sb{0}(0)}{g\sb{2}(1) g\sb{2}(0) }, \\
\end{split}
\end{equation}
and, 
\begin{equation}
\label{eq: twoidentity}
\begin{split}
\ \tfrac{G\sb{0}(1) }{ G\sb{2}(1)} -\eta\sb{G}\, \tfrac{G\sb{0}(0) }{ G\sb{2}(0)} 
&= \tfrac{H\sb{1} G\sb{2}(0) -H\sb{2} G\sb{0}(0) }{ G\sb{2}(1) G\sb{2}(0) } 
 +(1 -\eta\sb{G}) \tfrac{G\sb{0}(0)}{G\sb{2}(0)} \\
&\le \tfrac{H\sb{1} G\sb{2}(0) -H\sb{2} G\sb{0}(0) }{ G\sb{2}(1) G\sb{2}(0) } 
 +\tfrac{G\sb{0}\sp{2}(0)}{2 G\sb{2}\sp{2}(0)} \\
&=\tfrac{ G\sb{1}(1) G\sb{2}(0) -G\sb{1}(0) G\sb{2}(1)}{G\sb{2}(1) G\sb{2}(0)}
 +\tfrac{ (G\sb{1}(0) -G\sb{2}(0))\sp{2}}{2 G\sb{2}\sp{2}(0)}, \\
-\eta\sb{g}\, \tfrac{ g\sb{0}(1) }{ g\sb{2}(1) } +\tfrac{ g\sb{0}(0) }{ g\sb{2}(0) }
&= -\tfrac{h\sb{1} g\sb{2}(0) -h\sb{2} g\sb{0}(0)}{g\sb{2}(1) g\sb{2}(0)}
 +(1 -\eta\sb{g} ) \tfrac{g\sb{0}(0)}{g\sb{2}(0)} \\
&\le \tfrac{h\sb{2} g\sb{0}(0) -h\sb{1} g\sb{2}(0)}{g\sb{2}(1) g\sb{2}(0)}
 +\tfrac{g\sb{0}\sp{2}(0)}{2 g\sb{2}\sp{2}(0)} \\
&=\tfrac{g\sb{1}(0) g\sb{2}(1) -g\sb{1}(1) g\sb{2}(0)}{g\sb{2}(1) g\sb{2}(0)} 
 +\tfrac{ (g\sb{1}(0) -g\sb{2}(0))\sp{2}}{2 g\sb{2}\sp{2}(0)}, \\
\end{split}
\end{equation}
for the numerator and the denominator of ${\mathbf C}(s)$ in \eqref{eq: BC12defi}, 
respectively, recalling \eqref{eq: only10forl} with the first and the third 
inequalities in \eqref{eq: usednow} for the last but one step.

\subsection{Computation with details}. 
We recall the definitions of $G\sb{j}(l)$ and $g\sb{j}(l)$ for $j =1$, $2$ and 
$l =0$, $1$ in \eqref{eq: 2eqns} with \eqref{eq: ggGG1210}, getting 
\begin{equation*}
\label{eq: G2lwdA}
\begin{split}
G\sb{1}(1) &=G\sb{1}( {\hat\theta}\sb{1}, {\hat r}\sb{1}; k)  
=(W\sb{1} -\tfrac{1}{2})\sp{2} +{\hat r}\sb{1}\sp{2} 
 -2(W\sb{1} -\tfrac{1}{2}) {\hat r}\sb{1} 
  \cos({\hat\theta}\sb{1} -\tfrac{k\pi}{2\sp{K}} ), \\
G\sb{2}(0) &=G\sb{2}( {\hat\theta}\sb{0}, {\hat r}\sb{0}; k)   
=(W\sb{2} -\tfrac{1}{2})\sp{2} +{\hat r}\sb{0}\sp{2} 
 -2(W\sb{2} -\tfrac{1}{2}) {\hat r}\sb{0} 
  \cos({\hat\theta}\sb{0} -\tfrac{k\pi}{2\sp{K}} ). \\
\end{split}
\end{equation*}
Also, we recall the definitions of $W\sb{1} = 3R +R\sp{1/4} +\tfrac{1}{2}$ 
and $W\sb{2} =3R +\tfrac{1}{2}$ after \eqref{eq: prenabla}. The last 
expressions are simplified to 
\begin{equation}
\label{eq: G2lwdB}
\begin{split}
G\sb{1}(1) &= 9R\sp{2} +6 R\sp{5/4} +R\sp{1/2} +{\hat r}\sb{1}\sp{2} 
 - ( 6 R +2 R\sp{1/4} ){\hat r}\sb{1}\cos({\hat\theta}\sb{1} 
  -\tfrac{k \pi}{2\sp{K}} ) , \\
G\sb{2}(0) &= 9R\sp{2} +{\hat r}\sb{0}\sp{2} 
 -6 R{\hat r}\sb{0}\cos({\hat\theta}\sb{0} -\tfrac{k \pi}{2\sp{K}} ). \\
\end{split} 
\end{equation}
Similarly, we acquire 
\begin{equation}
\label{eq: G2lwd0}
\begin{split}
G\sb{1}(0) & =9R\sp{2} +6 R\sp{5/4} +R\sp{1/2} +{\hat r}\sb{0}\sp{2} 
 - ( 6 R +2 R\sp{1/4} ){\hat r}\sb{0} \cos({\hat\theta}\sb{0} 
  -\tfrac{k \pi}{2\sp{K}} ) , \\
G\sb{2}(1) &= 9R\sp{2} +{\hat r}\sb{1}\sp{2} 
 -6 R{\hat r}\sb{1} \cos({\hat\theta}\sb{1} -\tfrac{k \pi}{2\sp{K} } ) , \\
\end{split} 
\end{equation}
and 
\begin{equation}
\label{eq: g2lwd0}
\begin{split}
g\sb{1}(1) &= 9R\sp{2} +6 R\sp{5/4} +R\sp{1/2} +{\check r}\sb{1}\sp{2} 
 - ( 6 R +2 R\sp{1/4} ){\check r}\sb{1}\cos( {\check\theta}\sb{1} 
  -\tfrac{k \pi}{2\sp{K}} ) , \\
g\sb{2}(0) &= 9R\sp{2} +{\check r}\sb{0}\sp{2} 
 -6 R{\check r}\sb{0}\cos({\check\theta}\sb{0} -\tfrac{k \pi}{2\sp{K}} ), \\
g\sb{1}(0) & =9R\sp{2} +6 R\sp{5/4} +R\sp{1/2} +{\check r}\sb{0}\sp{2} 
 - ( 6 R +2 R\sp{1/4} ){\check r}\sb{0}\cos({\check\theta}\sb{0} 
  -\tfrac{k \pi}{2\sp{K}} ) , \\
g\sb{2}(1) &= 9R\sp{2} +{\check r}\sb{1}\sp{2} 
 -6 R{\check r}\sb{1}\cos({\check\theta}\sb{1} -\tfrac{k \pi}{2\sp{K}} ). \\
\end{split} 
\end{equation}

We recall the definition of $s\sb{0} =\sigma\sb{0} +i\,t\sb{0}$ with $\sigma\sb{0}
=\tfrac{5}{4}$ and $T -\tfrac{3}{2} <t\sb{0}\le T +\tfrac{3}{2}$ in the statement 
of Lemma \ref{lem: BC12upd}. It follows from \eqref{eq: eight} with 
\eqref{eq: calculatelater} that 
\begin{equation}
\label{eq: twotildeR}
\begin{split}
&\hskip 1.2true cm  
\sqrt{ \tfrac{81}{16} +(T -\tfrac{3}{2})\sp{2} }  
 =\sqrt{\tfrac{ 9\sp{2}}{4\sp{2}} +T\sp{2} -3 T +\tfrac{9}{4} } \\
&< T\,\sqrt{1 -\tfrac{3}{T} +\tfrac{117}{16 T\sp{2} } } 
 <{\hat r}\sb{0} =| s\sb{0} +1 | =\sqrt{\tfrac{81}{16} +t\sb{0}\sp{2} } 
  \\
&\hskip 1true cm
\le \sqrt{ \tfrac{81}{16} +(T +\tfrac{3}{2})\sp{2} } 
 <T\,\sqrt{1 +\tfrac{3}{T} +\tfrac{117}{16 T\sp{2} } }. \\
\end{split}
\end{equation}
Similarly, we get, for $\tfrac{1}{2} <X <1$, using the above 
restrictions on $\sigma\sb{0}$ and $t\sb{0}$: 
\begin{equation}
\label{eq: twotildeR2}
\begin{split}
&\hskip 2true cm 
T\,\sqrt{1 -\tfrac{3}{T} +\tfrac{45}{16 T\sp{2} } } 
<{\check r}\sb{0} =| s\sb{0} -(1 -X) | <T\,\sqrt{1 +\tfrac{3}{T} 
 +\tfrac{61}{16 T\sp{2} } }, \\
\end{split}
\end{equation}
For the estimates of ${\hat r}\sb{1}$ 
and ${\check r}\sb{1}$, one may use the triangle inequalities $||x| -|y|| 
\le |x +y|$ and $|x +y| \le |x| +|y|$ for any $x$, $y\in{\mathbb C}$. One 
acquires
\begin{equation}
\label{eq: sigmalr}
\begin{split}
|s\sb{0} +1| -1 \le &{\hat r}\sb{1} =|s\sb{1} +1| 
 \le |s\sb{0} +1| +1 , \\
|s\sb{0} -(1 -X)| -(1 -X) \le &{\check r}\sb{1} =|s\sb{1} -(1 -X)| \\
 &\le |s\sb{0} -(1 -X)| +(1 -X), \\
\end{split}
\end{equation}
using $| s\sb{1} -s\sb{0}| =1$ and $0< 1 -X <\tfrac{1}{2}$. With 
\eqref{eq: twotildeR}, \eqref{eq: twotildeR2}, and \eqref{eq: sigmalr} 
and using  $1 -\tfrac{1}{2} x< \sqrt{1 +x} <1 +\tfrac{1}{2} x$
for $0< x< 1$, we acquire 
\begin{equation}
\label{eq: fourr}
\begin{split}
T -\tfrac{3}{2} +\tfrac{117}{32 T} < 
 &{\hat r}\sb{0} <T +\tfrac{3}{2} +\tfrac{117}{32 T} , \\
T -\tfrac{1}{2} +\tfrac{117}{32 T}< 
 &{\hat r}\sb{1} <T +\tfrac{5}{2} +\tfrac{117}{32 T}, \\
T -\tfrac{3}{2} +\tfrac{45}{32 T}< 
 &{\check r}\sb{0} <T +\tfrac{3}{2} +\tfrac{61}{32 T}, \\
T -\tfrac{1}{2} +\tfrac{45}{32 T}< 
 &{\check r}\sb{1} <T +\tfrac{5}{2} +\tfrac{61}{32 T} . \\
\end{split}
\end{equation}

\subsection
{Proof of the first inequality in \eqref{eq: nablaCsDs} }.
\label{eq: 34}
 
Now, we estimate the expressions on the right hand sides of the last 
inequalities. We recall $2 T-1<  R\le 2 T+1$ from Section \ref{sec: intro} 
and use the estimates for ${\hat r}\sb{l}$ and ${\check r}\sb{l}$ for 
$l =0$, $1$ in \eqref{eq: fourr}. We also use $-1\le \cos(\theta)\le 1$ 
for any $\theta\in{\mathbb R}$. From \eqref{eq: G2lwdB}, we exhibit 
the bound on $G\sb{1}(1)$ in every detail. First, we group the terms 
involving higher powers of $T$, getting that  
\begin{equation}
\label{eq: 3stuffs}
\begin{split}
G\sb{1}(1) 
&\le 9(2 T +1)\sp{2} +6(2 T +1)\sp{5/4} +(2 T +1)\sp{1/2} 
 +(T +\tfrac{5}{2} \\ 
 &\hskip 1true cm  +\tfrac{117}{32 T})\sp{2} +[ 6(2 T +1) +2 (2 T +1)\sp{1/4} ] 
  (T +\tfrac{5}{2} +\tfrac{117}{32 T}) \\
 &= 9(2 T +1)\sp{2} +(T +\tfrac{5}{2} +\tfrac{117}{32 T})\sp{2} 
  +6(2 T +1)(T +\tfrac{5}{2} +\tfrac{117}{32 T}) \\
 &+ 6(2 T +1)\sp{5/4} + 2(2 T +1)\sp{1/4} (T +\tfrac{5}{2} 
  +\tfrac{117}{32 T}) +(2 T +1)\sp{1/2}. \\
\end{split}
\end{equation}
For the first three terms in the last sum, we have
\begin{equation}
\label{eq: 3stuffs01}
\begin{split}
&9(2 T +1)\sp{2} +(T +\tfrac{5}{2} +\tfrac{117}{32 T})\sp{2} 
  +6(2 T +1)(T +\tfrac{5}{2} +\tfrac{117}{32 T}) \\
&=49\,{T}^{2}+77\,T +\tfrac{1303}{16} +\tfrac{1287}{32\,T}
  +\tfrac{13689}{1024\,{T}^{2}}. \\  
\end{split}
\end{equation}
For the last three terms of the last expression in \eqref{eq: 3stuffs}, we 
use $(1 +x)\sp{\alpha} <1 +\alpha x$ for $0 <x <1$ and $0< \alpha<1$. We have 
\begin{equation}
\label{eq: 3stuffs02}
\begin{split}
&6(2 T +1)\sp{5/4} + 2(2 T +1)\sp{1/4} (T +\tfrac{5}{2} 
  +\tfrac{117}{32 T}) +(2 T +1)\sp{1/2} \\
&\le 6\times 2\sp{1/4} (2 T +1) T\sp{1/4} (1 +\tfrac{1}{8T})
 +2\sp{5/4} T\sp{1/4} ( 1+\tfrac{1}{8T})\\
 &\hskip 1true cm  \times(T +\tfrac{5}{2}+\tfrac{117}{32 T}) 
  +2\sp{1/2} T\sp{1/2} (1 +\tfrac{1}{4T}). \\  
\end{split}
\end{equation}
Expanding the last expressions in \eqref{eq: 3stuffs01} and \eqref{eq: 3stuffs02}
by a computer algebra package and putting back to \eqref{eq: 3stuffs}, we obtain
\begin{equation}
\label{eq: 3stuffs1}
\begin{split}
G\sb{1}(1) 
&\le 49\,{T}^{2} +7\times {2}^{{5}/{4}}\,{T}^{{5}/{4}} +77\,T +\sqrt{2}\,\sqrt{T}   
 +\tfrac{51\times{2}^{{1}/{4}}\,{T}^{{1}/{4}}}{8} \\
 &\quad +\tfrac{1303}{16} +\tfrac{{2}^{{1}/{2}}}{4\,\sqrt{T}} 
  +\tfrac{139\times{2}^{{1}/{4}}}{32\, {T}^{{3}/{4}}}
  +\tfrac{1287}{32\,T}+\tfrac{117\times{2}^{{1}/{4}}}{256\,{T}^{{7}/{4}}} 
  +\tfrac{13689}{1024\,{T}^{2}} \\
 &\le 49\,{T}^{2} +7\times {2}^{{5}/{4}}\,{T}^{{5}/{4}} +77\,T +\sqrt{2}\,\sqrt{T}   
 +15.228\,{T}^{{1}/{4}}. \\
\end{split}
\end{equation}
Here, we have only kept the first four terms with the higher powers of $T$ by 
noting that 
\begin{equation*}
\tfrac{51\times{2}^{{1}/{4}}\,{T}^{{1}/{4}}}{8}
 +\tfrac{1303}{16} +\tfrac{{2}^{{1}/{2}}}{4\,\sqrt{T}} 
  +\tfrac{139\times{2}^{{1}/{4}}}{32\, {T}^{{3}/{4}}}
  +\tfrac{1287}{32\,T}+\tfrac{117\times{2}^{{1}/{4}}}{256\,{T}^{{7}/{4}}} 
  +\tfrac{13689}{1024\,{T}^{2}} <15.228 T\sp{1/4},
\end{equation*}
as $T \ge T\sb{0}$, which is sufficient for our estimate in \eqref{eq: tonext0}
below. 

We do not exhibit the details in the following computation, since they are 
similar to ones which lead us to \eqref{eq: 3stuffs1}. From \eqref{eq: G2lwdB} 
and \eqref{eq: G2lwd0}, we obtain the bounds on both ends for $G\sb{2}(0)$ 
and $G\sb{1}(0)$, but only a lower bound for $G\sb{2}(1)$ is needed 
for \eqref{eq: tonext0} below. That is,
\begin{equation}
\label{eq: 3stuffs2}
\begin{split}
&\hskip 0.7true cm
49\,{T}^{2} +21\,T +53.437 < G\sb{2}(0) < 49\,{T}^{2} +63\,T +71.438,  \\
&\hskip 0.5true cm 49\,{T}^{2} +7\times {2}^{{5}/{4}}\,{T}^{{5}/{4}} 
 +21\,T+\sqrt{2}\,\sqrt{T} +5.691\,{T}^{{1}/{4}} < \\
& G\sb{1}(0) \le 49\,{T}^{2} +7\times {2}^{{5}/{4}}\,{T}^{{5}/{4}} 
 +63\,T+\sqrt{2}\,\sqrt{T} +12.842\,{T}^{{1}/{4}},  \\
&\hskip 2.5true cm G\sb{2}(1) >49\,{T}^{2} +35\,T +57.437; \\
&G\sb{1}(0) -G\sb{2}(0) <7\times {2}^{{5}/{4}}\,{T}^{{5}/{4}} 
 +42\,T+\sqrt{2}\,\sqrt{T} +12.842\,{T}^{{1}/{4}}. \\
\end{split}
\end{equation}

Recalling \eqref{eq: twoidentity} with \eqref{eq: 3stuffs1} and 
\eqref{eq: 3stuffs2}, we obtain
\begin{equation}
\label{eq: tonext0}
\begin{split}
&\tfrac{ G\sb{0}(1) }{ G\sb{2}(1) } -\eta\sb{G} \tfrac{ G\sb{0}(0) }{ G\sb{2}(0) }
 <\tfrac{(49\,{T}^{2}+7\,{2}^{{5}/{4}}\,{T}^{{5}/{4}}
 +\sqrt{2}\,\sqrt{T}+77.0001)(49\,{T}^{2}+63\,T+71.438)} 
 {\left( 49\,{T}^{2}+21\,T+54.437\right) \,
 \left( 49\,{T}^{2}+35\,T+54.437\right)} \\
&\hskip 1true cm
 -\tfrac{(49\,{T}^{2}+7\,{2}^{{5}/{4}}\,{T}^{{5}/{4}}
 +\sqrt{2}\,\sqrt{T}+22)(49\,{T}^{2}+35\,T+58)}
 {\left( 49\,{T}^{2}+21\,T+55\right) \,
 \left( 49\,{T}^{2}+35\,T+55)\right)} \\
&\hskip 2true cm  
 +\tfrac{(7\times{2}^{{5}/{4}}\,{T}^{{5}/{4}}+42\,T
 +\sqrt{2}\,\sqrt{T}+12.842\,{T}^{\frac{1}{4}})\sp{2}} 
 {2\left( 49\,{T}^{2}+21\,T+54.437\right)\sp{2}} \\
&\hskip 1true cm <\tfrac{1372\,{T}^{3}+ 937.5\,{T}^{{9}/{4}}}
 {2401\,{T}^{4}+2744.0\,{T}^{3}} 
 +\tfrac{49\,{2}^{{5}/{2}}\,{T}^{{5}/{2}}+1400}
 {4802.0\,{T}^{4}+4116.0\,{T}^{3}} <\tfrac{1}{2 T} .\\
\end{split}
\end{equation}

Similarly, from \eqref{eq: g2lwd0} we acquire the following estimates 
\begin{equation}
\label{eq: 3stuffsg}
\begin{split}
&g\sb{1}(0) < 49\,{T}^{2} +7\times {2}^{{5}/{4}}\,{T}^{{5}/{4}} 
 +63\,T+\sqrt{2}\,\sqrt{T} +12.9\,{T}^{{1}/{4}},  \\
&\hskip 0.6true cm 
 49\,{T}^{2} +35\,T +24.1 <g\sb{2}(1) <49\,{T}^{2} +77\,T +57.0 \\
&g\sb{1}(1) > 49\,{T}^{2} +7\times {2}^{{5}/{4}}\,{T}^{{5}/{4}} 
 +35\,T +\sqrt{2}\,\sqrt{T} +7.9\,{T}^{{1}/{4}}, \\
&\hskip 2.2true cm
 g\sb{2}(0) > 49\,{T}^{2} +21\,T +20.1;  \\
&g\sb{1}(0) -g\sb{2}(0) <7\times {2}^{{5}/{4}}\,{T}^{{5}/{4}} 
 +42\,T+\sqrt{2}\,\sqrt{T} +12.9\,{T}^{{1}/{4}}, \\
\end{split}
\end{equation}
which are all needed this time for \eqref{eq: tonext1} below.

It follows namely from \eqref{eq: twoidentity} with \eqref{eq: 3stuffs} that 
\begin{equation}
\label{eq: tonext1}
\begin{split}
&\eta\sb{g} \tfrac{ g\sb{0}(1) }{ g\sb{2}(1) } 
 -\tfrac{ g\sb{0}(0) }{ g\sb{2}(0) }
< \tfrac{(49\,{T}^{2}+7\,{2}^{{5}/{4}}\,{T}^{{5}/{4}}
 +63 T +\sqrt{2}\,\sqrt{T}+12.9)(49\,{T}^{2}+77\,T+57.0)} 
 {\left( 49\,{T}^{2}+77\,T+57.0\right) \,
 \left( 49\,{T}^{2}+21\,T+20.1\right)} \\
&\hskip 1true cm
  -\tfrac{(49\,{T}^{2}+7\,{2}^{{5}/{4}}\,{T}^{{5}/{4}}
 +35 T +\sqrt{2}\,\sqrt{T}+7.9)(49\,{T}^{2}+21\,T+20.1)}
 {\left( 49\,{T}^{2}+77\,T+57.0\right) \,
 \left( 49\,{T}^{2}+21\,T+20.1\right)} \\
&\hskip 2true cm 
 +\tfrac{(7\times{2}^{{5}/{4}}\,{T}^{{5}/{4}}+42\,T
 +\sqrt{2}\,\sqrt{T}+12.9\,{T}^{\frac{1}{4}})\sp{2}} 
 {2\left( 49\,{T}^{2}+21\,T+20.1\right)\sp{2}} \\
&\hskip 1true cm <\tfrac{4116\,{T}^{3}+937.5\,{T}^{{9}/{4}}}
{2401\,T\sp{4} +2744.0} 
 +\tfrac{49\,{2}^{{5}/{2}}\,{T}^{{5}/{2}}+1400.0}
 {4802.0\,{T}^{4}+4116.0\,{T}^{3}} <\tfrac{2}{T} .\\
\end{split}
\end{equation}

From \eqref{eq: lowerupper} and \eqref{eq: twoidentity} with \eqref{eq: tonext0} 
and \eqref{eq: tonext1}, we get, 
\begin{equation}
\label{eq: l7laststep}
\log\bigl| \tfrac{{\mathbf C}(s)}{{\mathbf C}(s\sb{0})} \bigr|
\le \tfrac{5 q}{4 T} <\tfrac{11 \log T}{T\sp{1/4} },  
\end{equation}
using the value of $q =\tfrac{(R+10)\log(R/2) +4\log \Omega}{4 {\grave\gamma} 
R\sp{1/4}}$ after \eqref{eq: prenabla} in Lemma \ref{lem: forLh} with 
$\omega=\xi(1/2)$, $\Omega =\tfrac{47.545}{\xi(1/2)}$ and $R \le 2 T +1$. 
The proof of  the first estimate in \eqref{eq: nablaCsDs} follows now directly 
from \eqref{eq: l7laststep}.

\subsection{Proof of the second inequality}.
To prove the second estimate in \eqref{eq: nablaCsDs}, one may use the same lines 
as for the first estimate in \eqref{eq: nablaCsDs}; but, a much easier way to prove 
the second estimate in \eqref{eq: nablaCsDs} is to have a recourse to the double 
symmetry property of $\nabla(s)$, as follows.

For the numerator function $\nabla[ s -i(Y -Y\sb{1}) ]$ of ${\mathbf D}(s)$, we 
note that $\nabla[s -i( Y -Y\sb{1} ) ] =\nabla( t- Y +Y\sb{1} - i \sigma)$ after 
we turn it by the angle of $\tfrac{3\pi}{2}$ counter-clockwisely. We notice that 
$\nabla( t - Y +Y\sb{1} - i\sigma ) =\nabla( 2 -u + v + i w)$, with $\sigma$, $t$, 
and $Y\sb{1} - Y$ replaced by $w$, $v$, and $2- u$;  also, $\nabla(2 -u + v +i w) 
=\nabla(2 -X +s)$ with $u= 2 -Y\sb{1} +Y =2 -\eta$ and $v + iw =s$. We remark 
that $0<\eta<\tfrac{1}{4}$, which is implied by (3.1) in \cite{Dh7}?, satisfies 
the condition for $X$ in the estimate for ${\mathbf C}(s)$. Hence, the same 
upper bound of the function ${\mathbf C}(s)$ is also an upper bound of 
the function ${\mathbf D}(s)$ on the same circle. Thus, the proof of the second 
bound in \eqref{eq: nablaCsDs} is achieved.
\qed


\vfill

\end{document}